\documentclass{article}

\usepackage[english]{babel}

\usepackage[letterpaper,top=2cm,bottom=3cm,left=3cm,right=3cm,marginparwidth=1.75cm]{geometry}

\usepackage{amsmath}
\usepackage{amssymb}
\usepackage{graphicx}
\usepackage{color}
\usepackage{amsfonts}
\usepackage{amscd}
\usepackage{mathrsfs}
\usepackage{bm}
\usepackage{enumerate}
\usepackage{algorithmic, algorithm}
\usepackage[colorlinks=true, allcolors=blue]{hyperref}

\newcommand{\bx}{\bm{x}}

\newcommand{\bn}{\mathbf{n}}
\newcommand{\by}{\bm{y}}

\newcommand{\p}{\partial}

\newcommand{\mathd}{\mathrm{d}}

\numberwithin{equation}{section}

\title{On the implementation of linear finite element method for nonlocal diffusion model over 2D domain}
\author{Zuoqiang Shi\thanks{Yau Mathematical Sciences Center, Tsinghua University, Beijing, China,
100084. \textit{Email: zqshi@tsinghua.edu.cn.}}}
\date{}

\begin{document}
\maketitle

\begin{abstract}
We propose an implementation of linear finite element method for nonlocal diffusion problem in 2D space. In the implementation, we reduce the integral from 4D to 2D which would simplify the computation significantly.
\end{abstract}

\section{Introduction}

For nonlocal diffusion model, finite element discretization usually has many good properties, such as asymptotic preserving, good convergence rate. However, direction implementation of finite element method for nonlocal diffusion model is difficult which typically requires the computation of integrals over $2d$ dimensional space \cite{cookbook-FEM}, $d$ is the dimension of the computational domain. In this note, we introduce an implementation of linear finite element method for nonlocal diffusion problem in 2D space. We reduce the integral from 4D to 2D which would simplify the computation significantly.

\section{Nonlocal model for Poisson equation}
We consider the Poisson equation with Neumann boundary condition.
\begin{align}
\label{eq:poisson}
-\Delta u(\bx)+u(\bx)&=f(\bx),\quad \bx\in \Omega\subset \mathbb{R}^2,\\
\frac{\p u}{\p \bn}(\bx)&=g(\bx),\quad \bx\in \p\Omega. \nonumber
\end{align}
A nonlocal counterpart of Poisson equation is given as follows 
\begin{eqnarray}
\label{eq:nonlocal_diffusion}
\frac{1}{\delta^2}\int_\Omega R_\delta(\bx,\by)(u(\bx)-u(\by))\mathd \by+\int_\Omega \bar{R}_\delta(\bx,\by)u(\by)\mathd \by=\int_\Omega \bar{R}_\delta(\bx,\by)f(\by)\mathd \by+2\int_{\p \Omega} \bar{R}_\delta(\bx,\by)g(\by)\mathd \tau_{\by},\quad
\end{eqnarray}

\begin{align}
\label{eq:kernel}
  R_\delta(\bx,\by)= C_\delta R\left(\frac{\|\bx-\by\|^2}{4\delta^2}\right),\quad 
\bar{R}_\delta(\bx,\by)= C_\delta\bar{R}\left(\frac{\|\bx-\by\|^2}{4\delta^2}\right){,}
\end{align}
where 
\begin{align}
\label{eq:kernelbar}
\bar{R}(r)=\int_{r}^{+\infty}R(s)\mathd s ,\quad
\bar{\bar{R}}(r)=\int_{r}^{+\infty}\bar{R}(s)\mathd s 
\end{align}
The constant $C_\delta=\alpha_2\delta^{-2}$ in \eqref{eq:kernel}
is a normalization factor so that
\begin{align}
\label{eq:normal-barR}
\int_{\mathbb{R}^2}\bar{R}_\delta(\bx,\by)\mathd \by=2\pi \alpha_2 \int_0^1 \bar{R}(r^2)r\mathd r=1.
\end{align}
We require that $R$ is compactly supported and piecewise polynomial, i.e.
\begin{itemize}
\item $R(r)=0, r\ge 1$ and $R$ is polynomial in $[0,1]$.
\end{itemize}

To get discretization, first we decompose the computational domain $\Omega$ to triangular mesh
$$\Omega=\bigcup_{i=1}^N T_i.$$ 
On the triangular mesh, we construct linear basis $$\phi_{i,k},\; i=1,2,\cdots,N,\;k=1,2,3.$$
$\phi_{i,1},\phi_{i,2},\phi_{i,3}$ are linear basis over triangle $T_i$ and vanish outside.

Then, we can get the discretization of the nonlocal diffusion model
\begin{align*}
\frac{1}{\delta^2}\sum_{j=1}^N\sum_{l=1}^3 c_{j,l}\int_\Omega\phi_{i,k}(\bx)&\int_\Omega R_\delta(\bx,\by)(\phi_{j,l}(\bx)-\phi_{j,l}(\by))\mathd \by\mathd \bx+\sum_{j=1}^N\sum_{l=1}^3 c_{j,l}\int_\Omega \phi_{i,k}(\bx)\int_\Omega\bar{R}_\delta(\bx,\by)\phi_{j,l}(\by)\mathd \by\mathd \bx\\
&=\int_\Omega \phi_{i,k}(\bx)\int_\Omega \bar{R}_\delta(\bx,\by)f(\by)\mathd \by \mathd \bx+2\int_\Omega \phi_{i,k}(\bx)\int_{\p \Omega} \bar{R}_\delta(\bx,\by)g(\by)\mathd \tau_{\by}\mathd \bx,\quad
\end{align*}
The direct computation of the stiffness matrix needs 4D numerical integration which is intractable in practice. After some derivation, we can reduce the 4D numerical integration to 2D which will simplify the computation tremendously. 
\section{Computation of stiffness matrix}
First, we compute the coefficient in the first term.
\begin{align*}
&\int_\Omega\phi_{i,k}(\bx)\int_\Omega R_\delta(\bx,\by)(\phi_{j,l}(\bx)-\phi_{j,l}(\by))\mathd \by\mathd \bx\\
=&\frac{1}{2}\int_\Omega\int_\Omega R_\delta(\bx,\by)(\phi_{i,k}(\bx)-\phi_{i,k}(\by))(\phi_{j,l}(\bx)-\phi_{j,l}(\by))\mathd \bx\mathd \by
\end{align*}
For $i=j$, we have
\begin{align*}
&\int_\Omega\int_\Omega R_\delta(\bx,\by)(\phi_{i,k}(\bx)-\phi_{i,k}(\by))(\phi_{i,l}(\bx)-\phi_{i,l}(\by))\mathd \bx\mathd \by\\
=& \int_{T_i}\int_{T_i} R_\delta(\bx,\by)(\phi_{i,k}(\bx)-\phi_{i,k}(\by))(\phi_{i,l}(\bx)-\phi_{i,l}(\by))\mathd \bx\mathd \by+2\int_{\Omega\backslash T_i}\int_{T_i} R_\delta(\bx,\by)\phi_{i,k}(\bx)\phi_{i,l}(\bx)\mathd \bx\mathd \by\\
=& \int_{T_i}\int_{T_i} R_\delta(\bx,\by)\bm{a}_{i,k}\cdot(\bx-\by)\bm{a}_{i,l}\cdot(\bx-\by)\mathd \bx\mathd \by+2\int_{T_i}\left(\int_{\Omega\backslash T_i} R_\delta(\bx,\by)\mathd \by\right)\phi_{i,k}^2(\bx)\mathd \bx
\end{align*}
In above calculation, we use the fact that $\phi_{i,k}$ are linear function on $T_i$ and vanish outside.

Above integral can be simplified further using integration by parts,
\begin{align*}
&\int_{T_i}\int_{T_i} R_\delta(\bx,\by)\bm{a}_{i,k}\cdot(\bx-\by)\bm{a}_{i,l}\cdot(\bx-\by)\mathd \bx\mathd \by\\
=& -2\delta^2 \int_{T_i}\int_{T_i} \bm{a}_{i,k}\cdot\nabla_{\bx} \bar{R}_\delta(\bx,\by)\bm{a}_{i,l}\cdot(\bx-\by)\mathd \bx\mathd \by\\
=& -2\delta^2 \int_{T_i}\int_{\partial T_i} \bm{a}_{i,k}\cdot\bm{n}(\bx) \bar{R}_\delta(\bx,\by)\bm{a}_{i,l}\cdot(\bx-\by)\mathd S_{\bx}\mathd \by+2\delta^2 \int_{T_i}\int_{T_i} \bm{a}_{i,k}\cdot\bm{a}_{i,l} \bar{R}_\delta(\bx,\by)\mathd \bx\mathd \by\\
=& 4\delta^4\int_{\partial T_i} \bm{a}_{i,k}\cdot\bm{n}(\bx)  \left(\int_{T_i}\bm{a}_{i,l}\cdot\nabla_{\by}\bar{\bar{R}}_\delta(\bx,\by)\mathd \by\right)\mathd S_{\bx}+2\delta^2\bm{a}_{i,k}\cdot\bm{a}_{i,l}  \int_{T_i}\int_{T_i} \bar{R}_\delta(\bx,\by)\mathd \bx\mathd \by\\
=& 4\delta^4\int_{\partial T_i} \bm{a}_{i,k}\cdot\bm{n}(\bx)  \left(\int_{\partial T_i}\bm{a}_{i,l}\cdot\bm{n}(\by)\bar{\bar{R}}_\delta(\bx,\by)\mathd S_{\by}\right)\mathd S_{\bx}+2\delta^2\bm{a}_{i,k}\cdot\bm{a}_{i,l}  \int_{T_i}\int_{T_i} \bar{R}_\delta(\bx,\by)\mathd \bx\mathd \by
\end{align*}
For the case $i\neq j$, 
\begin{align*}
&\int_\Omega\int_\Omega R_\delta(\bx,\by)(\phi_{i,k}(\bx)-\phi_{i,k}(\by))(\phi_{j,l}(\bx)-\phi_{j,l}(\by))\mathd \bx\mathd \by\\
=& \int_{T_i\cup T_j}\int_{T_i\cup T_j} R_\delta(\bx,\by)(\phi_{i,k}(\bx)-\phi_{i,k}(\by))(\phi_{j,l}(\bx)-\phi_{j,l}(\by))\mathd \bx\mathd \by\\
=& -2\int_{T_i}\int_{T_j} R_\delta(\bx,\by)\phi_{i,k}(\by)\phi_{j,l}(\bx)\mathd \bx\mathd \by\\
=& 2\int_{T_i}\int_{T_j} R_\delta(\bx,\by)\bm{a}_{i,k}\cdot(\bx-\by)\phi_{j,l}(\bx)\mathd \bx\mathd \by-2\int_{T_i}\left(\int_{T_j} R_\delta(\bx,\by)\bar{\phi}_{i,k}(\bx)\phi_{j,l}(\bx)\mathd \bx\right))\mathd \by\\
=& -4\delta^2\int_{T_i}\int_{T_j}\bm{a}_{i,k}\cdot\nabla_{\bx}\bar{R}_\delta(\bx,\by)\phi_{j,l}(\bx)\mathd \bx\mathd \by-2\int_{T_j}\bar{\phi}_{i,k}(\bx)\phi_{j,l}(\bx)\left(\int_{T_i} R_\delta(\bx,\by)\mathd \by\right)\mathd \bx
\end{align*}
where $\bar{\phi}_{i,k}$ denotes the extension of ${\phi}_{i,k}$ over $\mathbb{R}^2$.

Similarly, the first integral can be simpified by integration by parts further.
\begin{align*}
&-4\delta^2\int_{T_i}\int_{T_j}\bm{a}_{i,k}\cdot\nabla_{\bx}\bar{R}_\delta(\bx,\by)\phi_{j,l}(\bx)\mathd \bx\mathd \by\\
=& -4\delta^2 \int_{T_i}\int_{\partial T_j} \bm{a}_{i,k}\cdot\bm{n}(\bx) \bar{R}_\delta(\bx,\by)\phi_{j,l}(\bx)\mathd S_{\bx}\mathd \by+4\delta^2 \int_{T_i}\int_{T_j} \bm{a}_{i,k}\cdot\bm{a}_{j,l} \bar{R}_\delta(\bx,\by)\mathd \bx\mathd \by\\
=& -4\delta^2\int_{\partial T_j} \bm{a}_{i,k}\cdot\bm{n}(\bx) \phi_{j,l}(\bx)\left(\int_{T_i} \bar{R}_\delta(\bx,\by)\mathd \by\right)\mathd S_{\bx}+4\delta^2 \bm{a}_{i,k}\cdot\bm{a}_{j,l}\int_{T_j}\left(\int_{T_i}  \bar{R}_\delta(\bx,\by)\mathd \by\right)\mathd \bx
\end{align*}

The coefficients corresponding to the zero order term, we have
\begin{align*}
& \int_\Omega \phi_{i,k}(\bx)\int_\Omega\bar{R}_\delta(\bx,\by)\phi_{j,l}(\by)\mathd \by\mathd \bx\\
=&\int_{T_i} \phi_{i,k}(\bx)\int_{T_j}\bar{R}_\delta(\bx,\by)\phi_{j,l}(\by)\mathd \by\mathd \bx\\
=& \int_{T_j}\int_{T_i} \bar{R}_\delta(\bx,\by)\bm{a}_{i,k}\cdot(\bx-\by)\mathd \bx\phi_{j,l}(\by)\mathd \by+\int_{T_j}\left(\int_{T_i} \bar{R}_\delta(\bx,\by)\mathd \bx\right)\bar{\phi}_{i,k}(\by)\phi_{j,l}(\by)\mathd \by\\
=& 2\delta^2\int_{T_i}\int_{T_j}\bm{a}_{i,k}\cdot\nabla_{\by}\bar{\bar{R}}_\delta(\bx,\by)\phi_{j,l}(\by)\mathd \by\mathd \bx+\int_{T_j}\left(\int_{T_i} \bar{R}_\delta(\bx,\by)\mathd \bx\right)\bar{\phi}_{i,k}(\by)\phi_{j,l}(\by)\mathd \by
\end{align*}
and
\begin{align*}
&2\delta^2\int_{T_i}\int_{T_j}\bm{a}_{i,k}\cdot\nabla_{\by}\bar{\bar{R}}_\delta(\bx,\by)\phi_{j,l}(\by)\mathd \by\mathd \bx\\
=& 2\delta^2 \int_{T_i}\int_{\partial T_j} \bm{a}_{i,k}\cdot\bm{n}(\by) \bar{\bar{R}}_\delta(\bx,\by)\phi_{j,l}(\by)\mathd S_{\by}\mathd \bx-2\delta^2 \int_{T_i}\int_{T_j} \bm{a}_{i,k}\cdot\bm{a}_{j,l} \bar{\bar{R}}_\delta(\bx,\by)\mathd \by\mathd \bx\\
=& 2\delta^2\int_{\partial T_j} \bm{a}_{i,k}\cdot\bm{n}(\by) \phi_{j,l}(\by)\left(\int_{T_i} \bar{\bar{R}}_\delta(\bx,\by)\mathd \bx\right)\mathd S_{\by}-2\delta^2 \bm{a}_{i,k}\cdot\bm{a}_{j,l}\int_{T_j}\left(\int_{T_i}  \bar{\bar{R}}_\delta(\bx,\by)\mathd \bx\right)\mathd \by
\end{align*}
The source term can be computed in the same way.
\begin{align*}
&\int_\Omega \phi_{i,k}(\bx)\int_\Omega \bar{R}_\delta(\bx,\by)f(\by)\mathd \by \mathd \bx\\
=&\sum_{j=1}^N\sum_{l=1}^3 f_{j,l}\int_{T_i} \phi_{i,k}(\bx)\int_{T_j} \bar{R}_\delta(\bx,\by)\phi_{j,l}(\by)\mathd \by \mathd \bx
\end{align*}
Finally, we turn to deal with the boundary term.
\begin{align*}
&\int_\Omega \phi_{i,k}(\bx)\int_{\p \Omega} \bar{R}_\delta(\bx,\by)g(\by)\mathd \tau_{\by}\mathd \bx\\
=& \int_{\p \Omega}g(\by) \int_{T_i} \phi_{i,k}(\bx) \bar{R}_\delta(\bx,\by)\mathd \bx \mathd S_{\by}\\
=& \int_{\p \Omega}g(\by) \int_{T_i} \bm{a}_{i,k}\cdot(\bx-\by) \bar{R}_\delta(\bx,\by)\mathd \bx \mathd S_{\by}+\int_{\p \Omega}g(\by)\phi_{i,k}(\by) \int_{T_i}   \bar{R}_\delta(\bx,\by)\mathd \bx \mathd S_{\by}\\
=& -2\delta^2\int_{\p \Omega}g(\by) \int_{T_i} \bm{a}_{i,k}\cdot\nabla_{\bx} \bar{\bar{R}}_\delta(\bx,\by)\mathd \bx \mathd S_{\by}+\int_{\p \Omega}g(\by)\phi_{i,k}(\by) \int_{T_i}   \bar{R}_\delta(\bx,\by)\mathd \bx \mathd S_{\by}\\
=& -2\delta^2\int_{\p \Omega} \int_{\partial T_i} \bm{a}_{i,k}\cdot\bm{n}(\bx) g(\by)\bar{\bar{R}}_\delta(\bx,\by)\mathd S_{\bx} \mathd S_{\by}+\int_{\p \Omega}g(\by)\phi_{i,k}(\by) \left(\int_{T_i}   \bar{R}_\delta(\bx,\by)\mathd \bx\right) \mathd S_{\by}\\
\end{align*}
Based on above calculations, it can be found that the 4D integral can be reduced to 2D if we can compute the integral of kernel functions over triangles in explicit way. Fortunately, 
$$\int_{T_i} R_\delta(\bx,\by)\mathd \by,\; \int_{T_i} \bar{R}_\delta(\bx,\by)\mathd \by,\; \int_{T_i} \bar{\bar{R}}_\delta(\bx,\by)\mathd \by$$ 
can be given in explicit form based on the assumption that $R$ is polynomial on $[0,1]$.

\section{Integration over triangle}

\begin{figure}
\center
\includegraphics[width=0.35\textwidth]{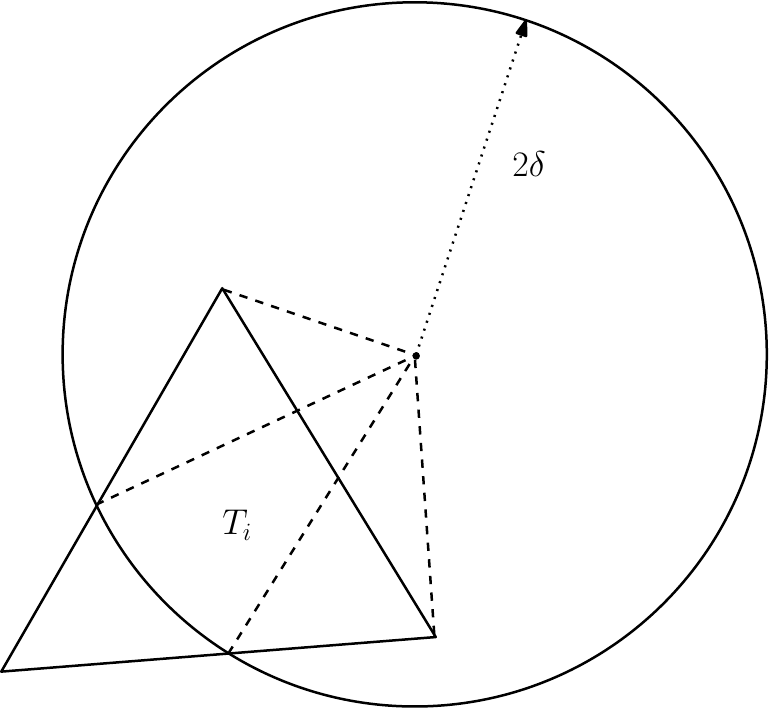}\hspace{15mm}
\includegraphics[width=0.33\textwidth]{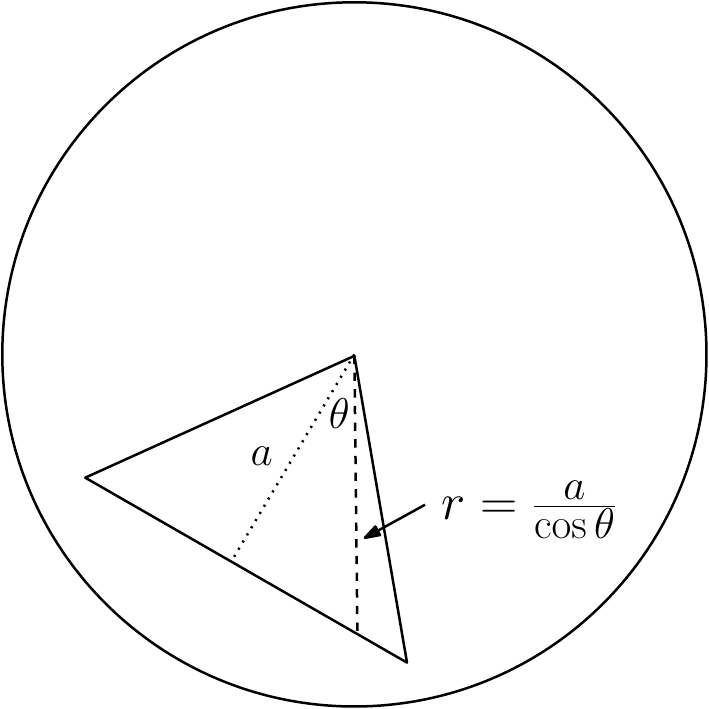}
\caption{Illustration of the integral over triangle. Left: decomposition of the integral domain; Right: integral over triangle.\label{fig:triangle}}
\end{figure}

As shown in the left figure of Fig. \ref{fig:triangle}, the integral over triangle can be decomposed to integral over circular sectors and triangles with center of the circle as one vertex. The integral over circular sectors is easy to calculate. The integral over triangles can be computed in radial coordinate as shown in the right figure of Fig. \ref{fig:triangle}. 

Since $R$ is polynomial on $[0,1]$, so 
\begin{align*}
\int_T R(\|\bx\|^2) \mathd \bx=\sum_k \int_T \|\bx\|^{2k} \mathd \bx
\end{align*}
For each term, we have
\begin{align*}
\int_T \|\bx\|^{2k} \mathd \bx=\int_{\theta_0}^{\theta_1} \int_0^{\frac{a}{\cos\theta}} r^{2k+1}\mathd r \mathd \theta= \frac{a^{2k+2}}{2k+2}\int_{\theta_0}^{\theta_1} (\sec\theta)^{2k+2}\mathd \theta
\end{align*}
On the other hand, using integration by part, we get
\begin{align*}
\int\sec^n \theta\mathd\theta=\frac{\sec^{n-2}\theta\tan \theta}{n-1}+\frac{n-2}{n-1}\int \sec^{n-2}\theta\mathd \theta,\quad n\neq 1.
\end{align*}
For $n=1,2$, it is easy to get
\begin{align*}
\int\sec \theta\mathd\theta=\ln |\sec \theta+\tan \theta|,\quad 
\int \sec^2\theta\mathd\theta=\tan\theta
\end{align*}
Finally, we can get explicit formula of $\int_{T_i} R_\delta(\bx,\by)\mathd \by$. For $\bar{R}$ and $\bar{\bar{R}}$, the explicit formula can be obtained in the similar way. 

\section{Conclusion}

In this note, we use nonlocal diffusion model with Neumann boundary condition as an example to introduce our method. This method can be easily extended to nonlocal diffusion model with Dirichlet boundary condition introduced in \cite{MP-dirichlet}. 
Next, we will try to extend the method to 3D case.

\bibliographystyle{alpha}
\bibliography{sample}

\end{document}